\documentstyle[12pt]{amsart}

\newcommand{\aut}{\operatorname{Aut}(G)}
\newcommand{\cobarasc}{{\cal Cobar} {\cal As^c}}
\newcommand{\cobarliec}{{\cal Cobar} {\cal Lie^c}}

\newcommand{\CHA}{$C_{\infty}$-algebra}

\newcommand{\CP}[1]{\Bbb{C}\Bbb{P}^{#1}}
\newcommand{\D}{{\cal D}}
\newcommand{\DKM}{De\-ligne-Knud\-sen-Mum\-ford\ }
\newcommand{\ec}{e_\nc}
\newcommand{\eg}{{\sl e.g.}}

\newcommand{\Ee}[1]{E^{#1}}

\newcommand{\End}[1]{{\cal E}@!n@!d\,_{#1}}
\newcommand{\F}[1]{F_{#1}}

\newcommand{\Hom}{\operatorname{Hom}}
\newcommand{\ie}{{\sl i.e.}}
\newcommand{\I}{\bf{I}}
\newcommand{\lb}{\lbrack}
\newcommand{\liec}{\operatorname{{\cal Lie^c}}}
\newcommand{\lie}{\operatorname{{\cal Lie}}}
\newcommand{\M}{\Mm}
\newcommand{\Mc}{\overline{\Mm}}

\newcommand{\Mm}{{\cal M}}

\newcommand{\nc}{{\Bbb{C}}}

\newcommand{\nz}{{\Bbb{Z}}}
\newcommand{\oo}{\circ}
\newcommand{\Oo}[1]{\OO({#1})}
\newcommand{\OO}{{ \cal O}}
\newcommand{\PSL}{\operatorname{PSL}(2,\nc)}

\newcommand{\rb}{\rbrack}

\newcommand{\SHA}{$A_{\infty}$-algebra}

\newtheorem{tth}{Theorem}[section]

\newtheorem{prop}[tth]{Proposition}
\newtheorem{crl}[tth]{Corollary}
\newtheorem{quest}[tth]{Question}

\theoremstyle{definition}

\newtheorem{df}{Definition}[section]
\newtheorem{ex}{Example}[section]

\theoremstyle{remark}

\newtheorem{rem}{Remark}
\newtheorem{ack}{Acknowledgment}

\begin{document}

\title[Moduli spaces and commutative homotopy algebras]
{Homology of moduli spaces of curves and commutative homotopy algebras}

\author
{Takashi Kimura}
\thanks{Research of the first author was supported in part by an NSF
Postdoctoral Research Fellowship}
\address
{Department of Mathematics, University of North Carolina, Chapel Hill, NC
27599-3250 }
\email{kimura@@math.unc.edu}

\author{Jim Stasheff}
\email{jds@@math.unc.edu}
\thanks{Research of the second author was supported in part by NSF grant
DMS-9206929}

\author
{Alexander A. Voronov}
\address
{Department of Mathematics, University of Pennsylvania, Philadelphia, PA
19104-6395}
\email{voronov@@math.upenn.edu}
\thanks{Research of the third author was supported in part by NSF grant
DMS-9402076}

\date{January 31, 1995}

\dedicatory{\smallskip
\hfill \protect\parbox{1.65in}{Now if it was dusk outside,\\
Inside it was dark.}
\smallskip\\
\hfill Robert Frost
\smallskip}

\maketitle

There have been a number of mathematical results recently identifying
algebras over certain operads \cite{bg:1, ksv, hs, konm, gk, g, gj,
grav}. See \cite {Ad, km}  for expository surveys of the basics of operad
theory. Before citing any of these results, let us mention some trivial
classical examples. Let ${\cal A}$ denote one of the three words:
``commutative'', ``associative'' and ``Lie''. In each of these cases,
consider the corresponding operad ${\cal O}(n) = \langle $ words in the free
${\cal A}$ algebra on $n$ generators having at most one occurance of  each
generator $ \rangle_k, \; n \ge 1$, $\langle \; \rangle_k$ meaning the linear
span over the ground field $k$. When ${\cal A}$ is commutative, associative or
Lie, we will denote the corresponding operad ${\cal O} $ by ${\cal C} $,
${\cal A} s $ and ${\cal L}ie $, respectively. Note that ${\cal C} (n) =
\langle x_1 \dots x_n \rangle_k = k$ and ${\cal A} s (n) = \langle
x_{\sigma(1)} \dots x_{\sigma(n)} \; | \; \sigma \in S_n \rangle_k=
k[S_n]$. The main feature of these three operads is that they describe
algebras of the corresponding types. More precisely, algebras over an operad
${\cal O}$, ${\cal O} = {\cal C}$, ${\cal As}$ or ${\cal Lie}$, (or simply,
${\cal O}$-algebras) are exactly ${\cal A}$ algebras.

A hint of another relation between
operads and algebras may be given by the formula:
\[
\bigoplus_{n \ge 1} ({\cal O} (n) \otimes V^{\otimes n})_{S_n} \text{
is the free ${\cal O}$-algebra generated by a vector space $V$.}
\]
Here $S_n$ in the subscript denotes coinvariants of the symmetric group, that
is to say, the quotient by the diagonal action of the symmetric group.

Several recent, less obvious examples were largely inspired by
the work of physicists, who came up with new and not so new kinds of
algebras, surprisingly appearing in certain models of quantum field
theory. An essential feature of these algebras was that they did not satisfy
the classical identities strictly, but rather satisfied them `up to homotopy'
in a strong sense. Mathematicians, in their attempt to describe the algebras
geometrically, found operads responsible for the types of algebras
produced by the physical theories, see Table~1.

\begingroup
\input{psfig}
\begin{figure}[h]
\vskip -.5truein\vbox{\hskip .75truein
\hbox{\psfig{figure=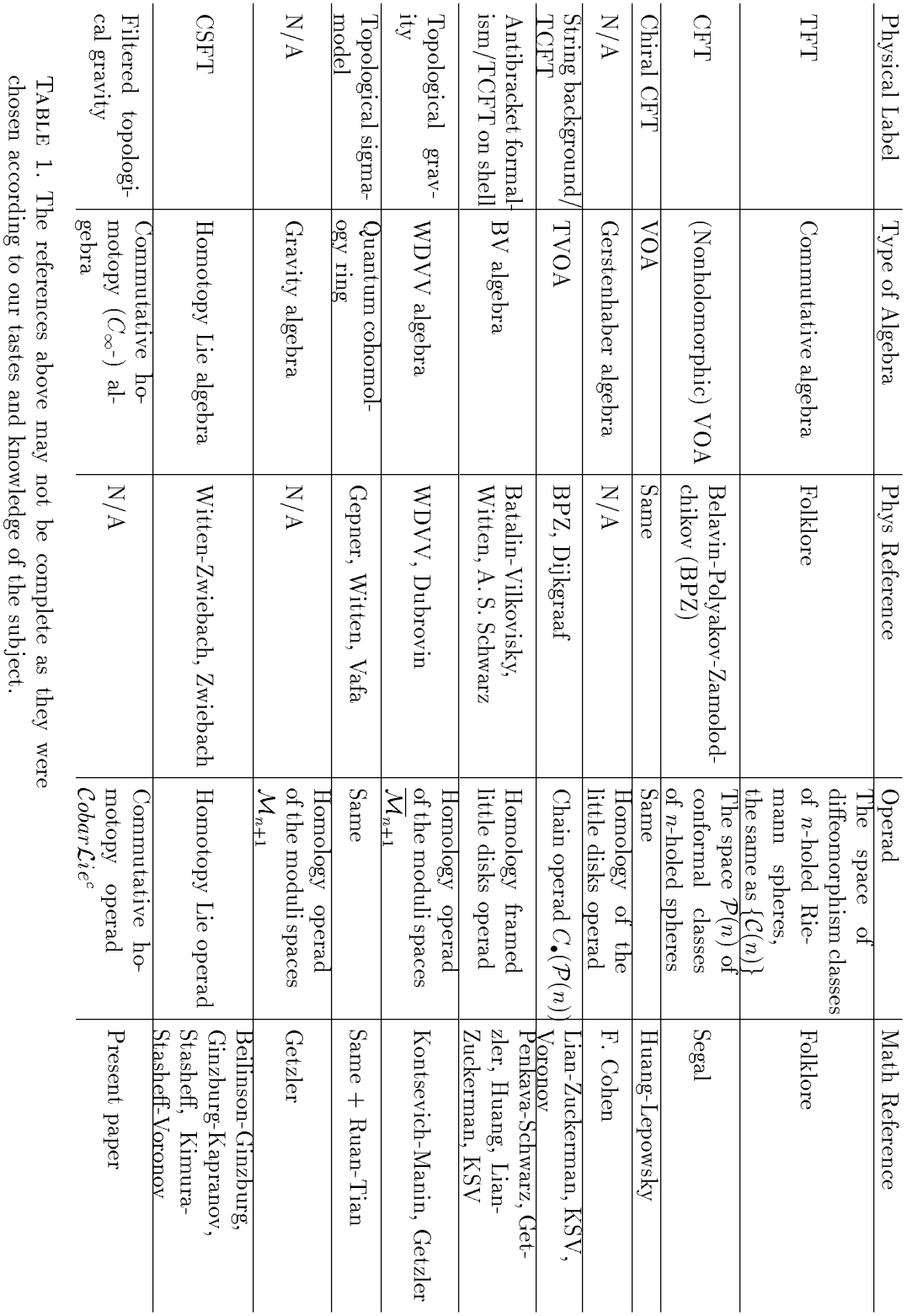,height=9in}}}
\label{table}
\end{figure}
\endgroup

This paper proposes a theory with a physics flavor, filtered topological
gravity, whose state space is endowed with the structure of a commutative
homotopy algebra or  $C_\infty$-algebra. The core of the paper is the
following theorem. Consider the operad of moduli spaces $\Mc(n)$, the \DKM
compactification of the moduli space of $(n+1)$-punctured Riemann
spheres. Let $C(\Mc(n))$ be the corresponding operad of singular chains. This
is an operad of complexes. Hence, an algebra over it is a complex $V$ with a
morphism of operads $C(\Mc(n)) \to \Hom (V^{\otimes n}, V)$.

\begin{tth}
\label{th:goal}
Let $V$ be an algebra over $C(\Mc(n))$ such that the structure morphism
$\mu:C(\Mc(n)) \,\to\, \Hom(V^{\otimes},V)$ vanishes on all elements in
$C_p(\Mc(n))$ for all $p > n-2$, then $V$ has the structure of a
$C_\infty$-algebra.
\end{tth}

This result may be regarded as a partial answer to a question of Lian
and Zuckerman \cite{lz} of lifting the BV algebra structure on BRST
cohomology to the cochain level, more precisely, providing the BV dot
product with higher homotopies. In fact, in \cite{k}, it was shown
that the BRST complex formed a homotopy associative algebra which
arose from an operad which is closely related to the commutative
homotopy associative operad above.

Another result, which is a byproduct of our work, is our description
of a spectral sequence  which
converges to $H^\bullet(\M_{g,n})$, whose $E_1$ term consists of tensor
products of the cohomology groups of compactified moduli spaces. Our result can
be stated
succinctly by using the generalization of the operad cobar
construction which includes graphs as well as trees, called the
Feynman transform, a notion due to Getzler and Kapranov
\cite{modular}.

\begin{tth}
There exists a spectral sequence converging to $H^\bullet (\M_{g,n})$ and
degenerating at the $E_2$ term, such that the term $E_1$ is the cobar
construction of the modular co-operad $H^\bullet (\Mc_{g,n})$.
\end{tth}

The proof of this theorem is based on purity of the Hodge structure on
$H^\bullet (\Mc_{g,n})$. We also present a ``dual'' version of this
theorem, where the spectral sequence converges to $H_\bullet
(\Mc_{g,n})$ and does not necessarily degenerate at $E^2$, unless
$g=0$.

\section{$C_\infty$ operad and $C_\infty$-algebras}

Let us begin with a definition of an operad.

\begin{df}
An {\sl operad  $\OO = \{\,\OO(n)\,\}_{n\geq 0}$ with unit} is a
collection of objects (topological spaces, complexes, etc.)  such that each
$\OO(n)$ has an action of $S_n$, the permutation group on $n$ elements ($S_0$
is contains only the identity), and a collection of operations for $n\geq 1$
and $1\leq i\leq n$, $\OO(n)\times\OO(n')\,\to\,\OO(n+n'-1)$ given by
$(f,f')\,\mapsto\,f\oo_i f'$ satisfying
\begin{enumerate}
\item if $f\in\OO(n),$ $f'\in\OO(n')$, and $f''\in\OO(n'')$ where $1\leq i <
j \leq n$, $n', n''\geq 0,$ and $n\geq 2$ then
\begin{equation}
\label{eq:nosign}
(f\oo_i f')\oo_{j+n'-1} f'' = (-1)^{|f'| |f''|} (f\oo_j f'')\oo_i f'
\end{equation}
where signs on the right hand side should be ignored if $\OO$ is not an
operad of (graded) vector spaces,
\item if $f\in\OO(n),$ $f'\in\OO(n')$, and $f''\in\OO(n'')$ where $n,n'\geq
1,$ $n''\geq 0$, and $i = 1, \ldots, n$ and $j = 1, \ldots, n'$ then
\begin{equation*}
(f\oo_i f')\oo_{i+j-1} f'' = f\oo_i (f'\oo_j f''),
\end{equation*}
\item the composition maps are equivariant under the action of the
permutation groups,
\item there exists an element $\I$ in $\OO(1)$ called the {\sl
unit} such that for all $f$ in $\OO(n)$ and $i=1,\ldots,n$,
\begin{equation*}
\I\oo_1 f = f = f\oo_i \I
\end{equation*}
\end{enumerate}
\end{df}

By iterating $k$ composition maps, one obtains the more common form of the
operad composition, $\gamma\,:\,\Oo{k} \times \Oo{n_1} \times\cdots\times
\Oo{n_k}\,\to\,\Oo{n_1 + \cdots + n_k}.$

\begin{ex} The {\sl endomorphism operad $\End{V}:=\{\End{V}(n)\}$ of a
differential graded $k$-vector space $(V,Q)$} is defined to be $\End{V}(n) :=
\Hom_k(V^{\otimes n},V)$ where $S_n$ acts naturally upon $V^{\otimes n}$ and
the composition maps are given by
\begin{multline}
 (f\oo_i f')(v_1\otimes \ldots \otimes v_{n+n'-1}) \\
 =\, \pm
f(v_1\otimes\ldots\otimes v_{i-1}\otimes f'(v_i\otimes\ldots v_{i+n'})\otimes
v_{i+n'+1}\otimes v_{n+n'-1})
\end{multline}
for all $f$ in $\End{V}(n)$, $f'$ in $\End{V}(n')$ and $i = 1,\ldots, n$
for all $n,n'$, and the unit $e$ is just the identity map from $V$ to itself.
The $\pm$ is the sign which is obtained by sliding $f'$ through $v_1, \ldots,
v_{i-1}$. Let us denote the component of $\End{V}(n)$ with degree $g$ by
$\End{V}^g(n)$.
\end{ex}

Algebraic structures on a differential graded vector space $V$ are often
parametrized by an operad through the following notion.

\begin{df}
Let $\OO$ be an operad. A differential graded $k$-vector space,
$(V,Q)$, is said to be an {\sl $\OO$-algebra} if there is a morphism of
operads $\OO\,\to\,\End{V}$.
\end{df}

Notice that given an operad of topological spaces $\OO$, the singular chains
on $\OO$, $C_\bullet(\OO):=\{C_\bullet(\OO(n))\}$,  naturally forms an operad
of differential graded vector spaces and, consequently, so do the homology
groups $H_\bullet(\OO):=\{\,H_\bullet(\OO(n))\,\}.$  Furthermore, if $(V,Q)$ is
an algebra over $C_\bullet(\OO)$, then $H_\bullet(V)$ is naturally an algebra
over $H_\bullet(\OO)$. However, there is more information in the
$C_\bullet(\OO)$-algebra, $(V,Q),$ than at the cohomology level. This
provides motivation for the study of algebraic structures up to homotopy, the
first of which we now recall.

\begin{df}
An {\sl $A_\infty$ (or strongly homotopy associative) algebra}, $V$, is a
complex $Q:V^g\,\to\,V^{g-1}$ endowed with a collection of $n$-ary (linear)
operations $\{\,m_n:V^{\otimes n}\,\to\,V\,\}_{n\geq 2}$ with $m_n$ having
degree $n-2$ satisfying
\begin{multline}
\label{eq:Aia}
Q (m_n(v_1,\ldots,v_n)) - (-1)^n \sum_{k=1}^n (-1)^{\epsilon(k)}
m_n(v_1,\ldots , Q v_k, \ldots, v_i) \\ = \sum_{r,s,k} (-1)^{k(s-1) + s n} \,
(m_r\oo_k m_s) (v_1,\ldots, v_n)
\end{multline}
where the summation runs over $r,s,k$ satisfying $r+s = n+1$, $ 1\leq k\leq r$,
$2\leq r < n$ for all $v_1,\ldots,v_n$ in $V$ and $n\geq 2$, and
$\epsilon(k)$ denotes the sign obtained by sliding $Q$ through $v_1,\ldots,
v_{k-1}$.
\end{df}

We note that $m_2$ is a multiplication, $m_3$ is an associating homotopy and
the  further $m_k$'s are ``higher associating homotopies''.

In the context of (strong) homotopy algebras, the simplest definition of a
\CHA,   first appeared in work of Kadeishvili \cite{kad1, kad2} and then in
that of Smirnov  \cite{smirn}  (both of whom called them commutative
$A_\infty$-algebras),  is:

\begin{df} A {\it \CHA}  is an \SHA\  $(A, \{ m_n\})$ such that each map
$m_n:A^{\otimes n}\to A$ is a Harrison cochain, i.e. $m_n$ vanishes on the
sum of all $(p,q)$-shuffles for $p+q=n$, the sign of the shuffle coming
from the grading of $A$ shifted by $1$.
\end{df}

The single object equivalent definition  (cf.\ \cite{SS, Q, jcm}) is:

\begin{df} A {\it \CHA} is a graded vector space $A$ together with a
codifferential on the free Lie coalgebra cogenerated by $sA$,
which is isomorphic to $A$ with the grading shifted by $1$.  (The
convention is that the shift is opposite to the degree of the differential;
since we will be working homologically at the operad level, $d$ is of
degree $-1$ while $s$ is of degree $+1$.)
\end{df}

This definition hints of  the ancient Koszul duality (or adjointness) between
commutative algebras and  Lie coalgebras (or Lie algebras and commutative
coalgebras) \cite{d,Q, jcm, SS}, which carries over to the operad level,
\cite{gk}.

Getzler and Jones \cite{gj} established the equivalence of the above
with  the operadic definition of a \CHA \ implicit in \cite{gk}.

\begin{df}  A {\it \CHA}  is an algebra over the operad $\cobarliec$.
\end{df}

The concept of a co-operad is defined dually to that of an operad \cite {gj},
so that if ${\cal K} = \{ {\cal K}(n) \}$ is a co-operad over a field $k$,
then the linear duals $Hom({\cal K}(n),k)$ form an operad.  Conversely, if
${\cal O} = \{{\cal O}(n) \}$ is an operad with each ${\cal O}(n)$ finite
dimensional  (or of finite type, i.e.  graded and finite dimensional in each
grading), then the linear duals $Hom({\cal O}(n), k)$ form a co-operad.  (A
co-operad  is ``an operad with the arrows reversed''.)

The operad $\lie = \{\lie (n)\}$ \cite{hs} with $\lie (n)$ is finite
dimensional for each n.  The co-operad $\liec$ has $\liec (n) = \Hom (\lie
(n), k).$

${\cal Cobar} $ \cite{gj} is a functor from co-operads to operads. (In
\cite{gk}, the use of co-operads is avoided by assuming finite type for
operads and defining ${\cal Cobar}$ only for the linear duals of operads.) All
we need to know is:

For any co-operad ${\cal K}$, ${\cal Cobar}\  {\cal K}$ is an operad with
pieces indexed by trees, constructed as products of various ${\cal K}(n)$'s
according to a prescription determined by the tree. Our convention is that
trees have vertices all of valence $ > 1$, e.g. the {\sl corolla} with $n$
leaves, denoted by $\delta_n$, and one root has just one vertex. The piece of
${\cal Cobar} \ {\cal K}$ indexed by this corolla is just ${\cal K}(n)$ with a
shift in grading.

\begingroup
\input{psfig}
\begin{figure}[h]
\centerline{\psfig{figure=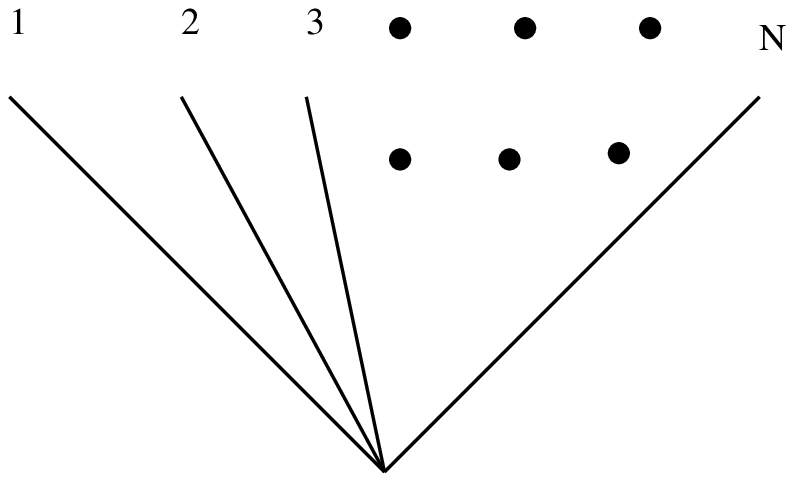,height=1in}}
\caption{An $N$-corolla, $\delta_N$}
\label{corolla}
\end{figure}
\endgroup

To make a comparison between \CHA s and \SHA s at the operad level, recall
that  the free Lie algebra ${\cal L}(x_1,\dots,x_n)$ can be realized as the
primitive subspace of $T(x_1,\dots,x_n)$ with  respect to the unshuffle
coproduct \cite{Ree}. Dually, the free Lie coalgebra ${\cal
L}^c(x_1,\dots,x_n)$ can then be identified with the  space of
indecomposables of the free associative coalgebra $T^c(x_1,\dots,x_n)$,
i.e. the {\it quotient\/} by the
image of the shuffle product. (See  \cite{WM} for definitions which do not rely
on finite type and duality.) Then ${\cal Lie}^c (n)$ is defined as the {\it
quotient\/} of the free Lie coalgebra ${\cal L}^c(x_1,\dots,x_n)$  by those
multilinear functions which vanish whenever two arguments are equal.

An \SHA\  is an algebra over the operad ${\cal Cobar} {\cal As}^c$ where
${\cal As}^c$
is the co-operad for associative coalgebras.  ${\cal As}^c (n)$ can be
realized as the quotient of the tensor coalgebra $T^c (x_1, \dots, x_n)$
by those multilinear functions which vanish whenever two arguments are equal.

Since $\cobarliec (n)$ has its corolla component equal to ${\cal Lie}^c (n)$
with a shift in grading, the structure map $m_n$ and its permutations are
given by  ${\cal Lie}^c (n) \to \Hom (A^{\otimes n}, A)$.  Interpretation of
$m_n$ as the structure map for an \SHA\  means pulling back the map
${\cal Lie}^c (n) \to \Hom (A^{\otimes n}, A)$ to $ T^c (n) \to \Hom(A^{\otimes
n}, A)$, which guarantees that $m_n$ vanishes on the image of the shuffle
product.

Strictly speaking, an algebra over ${\cal As}$, i.e. a morphism ${\cal As}\to
End (V)$, does not determine a unique associative multiplication on $V$, but
rather two of them, since ${\cal As}(2) = k\lb S_2\rb.$  We make the
obvious choice corresponding to the identity element of $S_2$ to
determine an associative multiplication on $V$. For ${\cal Cobar} {\cal As}^c$,
we must make choices for each $m_n$, but  we make the same choice,
corresponding to the identity element of $S_n$.  (This choice is implicit in
the result of \cite{gk} that an \SHA\  is the same as an algebra over
$\cobarasc.$) For $\cobarliec$, we just take the equivalence class (mod the
shuffle product) of the above choice for $\cobarasc$.

\section{Moduli Spaces of Punctured Spheres}

Let $\Mm_n$ be the moduli space of Riemann spheres with $n$ punctures. That
is, points in $\Mm_n$ consist of configurations of $n$ ordered points on
$\CP{1}$ with any two such configurations being identified if they are
related by a biholomorphic map. In other words,
$$\Mm_n := ((\CP{1})^n \setminus \Delta)/ \PSL$$
where $\Delta = \{(z_1,\dots , z_n) \linebreak[0] \in (\CP{1})^n \; |\; z_i =
z_j \; \text {\ for some \ } \; i \ne j\}$, the set of diagonals.
There is a compactification of $\Mm_n$ when $n\geq 3$ due to \DKM
\cite{dl,kap,keel,kn} which is the {\sl moduli space of stable genus 0 curves
with $n$ punctures,} $\Mc_n$.

Recall that a stable $n$ punctured complex curve of genus 0 is a connected
compact complex curve $C$ of genus 0 with $n$ punctures, such that it may
have ordinary double points away from the punctures, each irreducible
component of the curve $C$ is a projective line and the total number of
punctures and double points on each component of $C$ is at least 3. Both
$\Mm_n$ and $\Mc_n$ are smooth complex algebraic manifolds of complex
dimension $n-3$. The moduli space $\Mm_n$ of nonsingular curves is an open
submanifold in the projective manifold $\Mc_n$. The complement is a divisor,
formed by all degenerate curves.

Let $\Mc(1) := \{\ec\}$ and $\Mc(n) := \Mc_{n+1}$ for $n\geq 2$, then the
set $\Mc\, :=\, \{\,\Mc(n)\,\}$ is naturally an operad of algebraic varieties
where the element $\ec$ is defined to be a unit with respect to the operad
composition.\footnote{We include a unit in this manner for convenience.} The
permutation group on $n$ elements, $S_n$, acts on $\Mc(n)$ by reordering the
first n-punctures. The composition maps $\gamma_i\, :\,
\Mc(n)\,\times\,\Mc(n')\,\to\,\Mc(n+n'-1)$ for all $i = 1, \ldots, n$ and for
all $n,n'$ are defined by
$$(\Sigma,\Sigma')\,\mapsto\, \gamma_i(\Sigma,\Sigma')\, :=\, \Sigma \oo_i
\Sigma'$$
where $\Sigma\oo_i\Sigma'$ is obtained by attaching the $(n'+1)$st puncture
of $\Sigma'$ to the $i$th puncture of $\Sigma$ thereby creating a curve with
a new double point and the $(n+n'-1)$ punctures are ordered in the natural
way.

Each space $\Mc(n)$ is naturally stratified by smooth, connected locally
closed algebraic subvarieties. Each stratum of $\Mc(n)$ consists of those
points arising from $n$-punctured stable curves of a given topological type.
Since any stable curve can be obtained by attaching spheres together, the
combinatorics of this attaching process can be encoded in a tree.  Therefore,
each stratum of $\Mc(n)$ is naturally indexed by a tree (see Figure
\ref{strata}). This stratification gives rise to a filtration of $\Mc(n)$:
$$ \F{-1}(n)\, =\, \emptyset \subseteq \F{0}(n) \subseteq \cdots \subseteq
\F{n-2}(n)\,=\, \Mc(n)$$
where $\F{p}(n)$ is the disjoint union of all the strata of $\Mc(n)$ with
complex dimension less than or equal to $p$.  Furthermore, this filtration
is  invariant  under the action of the permutation group  and it behaves
nicely with respect to operad composition, {\sl i.e.} for all
$i=1,\ldots, n$ and for all $n,n',p,p'$, we have
$$ \gamma_i\,:\,\F{p}(n) \times \F{p'}(n')\,\to\,\F{p+p'}(n+n'-1).$$
The filtration gives rise to a spectral sequence associated to each $\Mc(n)$
which converges finitely to $H_\bullet(\Mc(n))$.  This spectral sequence is
known to degenerate at the $E^2$-term \cite{bg:1}. We show in
Section \ref{sec:filter} that any filtration which respects the operad
structure induces an operad structure on the $\Ee{r}$ term in its associated
spectral sequence for all $r\geq 0$. In particular, each $E^r$ term contains
a collection of suboperads. In our case, the only nontrivial suboperad is
the $q=0$ (``middle'') row of the $\Ee{1}$ term. The main result for our
purposes, due to Beilinson-Ginzburg \cite{bg:1}, cf.\ F.~Cohen
\cite{fc} and Schechtman-Varchenko \cite{sv}, is the following:

\begin{tth}
The ``middle'' row of the $\Ee{1}$ term of the spectral sequence associated
to the canonical filtration of $\Mc(n)$
$$
0 @>>> H_{n-2} (\F{n-2}(n),\F{n-3}(n)) @>>>	\dots @>>>
H_p(\F{p}(n),\F{p-1}(n)) @>>> \dots @>>> H_0(\F{0}(n)) @>>> 0, $$
is an operad isomorphic to $\cobarliec (n)$, the $C_\infty$ operad.
\end{tth}

\begin{sloppypar}
This identification is quite explicit.  Lefshetz duality gives an
isomorphism $H_p(F_p(n),F_{p-1}(n),\nc) \simeq H^p(F_p(n) \setminus
F_{p-1}(n))$ but $F_p(n)\setminus F_{p-1}(n)$ is the disjoint union of
those strata in $\Mc(n)$ with complex dimension $p$ each of which is
indexed by an $n$-tree with $n-2-p$ internal edges. If $T$ is a tree, then
let $S_T$ denote the strata associated to this tree, \eg\ see Figure
\ref{strata}. Each such tree $T$ is the iterated composition of corollas,
say $\delta_{n_1},\ldots, \delta_{n_k}$, and we have the isomorphism $S_T
\simeq S_{\delta_{n_1}}\times\cdots\times S_{\delta_{n_k}}.$ The stratum
associated to any corolla, $\delta_n,$ can be identified with $\M(n)$.
Therefore, $S_T \simeq \M(n_1)\times\cdots\times \M(n_k)$. Using the fact
that the cohomology of $\M(n)$ vanishes in degree above $n-2$, which
follows from a similar vanishing of homology of configuration spaces by
Arnold \cite{Ar}, we obtain the isomorphism $H^p(S_T) \simeq
H^{n_1-2}(\M(n_1)) \otimes \cdots \otimes H^{n_k-2}(\M(n_k))$. However,
$\lie(n)$ can be identified with $H_{n-2}(\M(n))$ (with a shift in degree)
\cite{fc} and therefore, $H^{n-2}(\M(n))$ can be identified with
$\liec(n)$. Finally, we obtain $H^p(S_T) \simeq \liec(n_1)\otimes \cdots
\otimes \liec(n_k)$ which is an element of $\cobarliec(n)$ associated to
the tree $T$ with the proper element in $\liec$ decorating the
corresponding vertex of $T$.
\end{sloppypar}

\begingroup
\input{psfig}
\begin{figure}[h]
\centerline{\psfig{figure=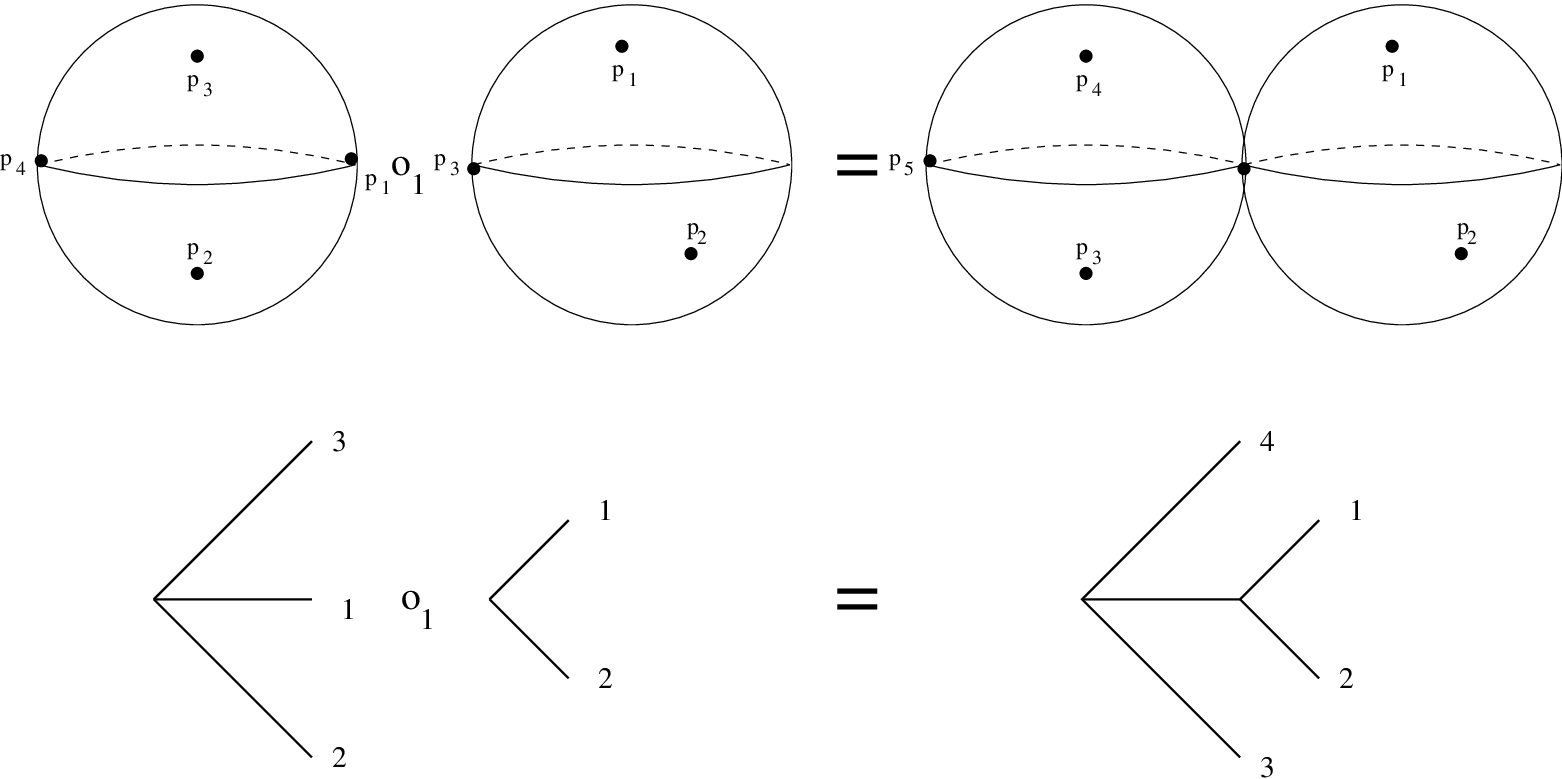,height=2in}}
\caption{The composition map $\oo_1:\,\Mc(3)\times\Mc(2)\,\to\, \Mc(4)$ and
the trees indexing the strata to which they belong}
\label{strata}
\end{figure}
\endgroup

\section{Moduli Spaces of Punctured Riemann Surfaces}

Here we make a generalization of the results of the previous section; we give
a complete description of the analogous spectral sequence in the case of moduli
spaces of Riemann surfaces of higher genera. In particular, we identify the
other rows of the $E^1$ term of the spectral sequence of the previous
section.

Let us restrict ourselves to those Riemann surfaces with genus $g$ and $n$
ordered punctures which have negative Euler characteristic,  \ie\ $g \ge 2$
or $g=1,
n \ge 1$ or $g = 0 , n \ge 3$. Let $\M_{g,n}$ be the moduli space of genus
$g$ Riemann surfaces, that  is, complex algebraic curves, with $n$ punctures
and let $\Mc_{g,n}$ be its compactification due to
Deligne-Knudsen-Mumford. This is a smooth, complete stack of dimension
$\dim_{\nc} \Mc_{g,n} = 3g - 3 + n$. The compactified moduli parameterize
isomorphism classes of stable curves, ones which have a finite number of
singularities, which are double points, and such that each irreducible
component of genus 1 has at least one puncture or double point and each
irreducible component of genus 0 has at least three punctures or double
points. These conditions insure that each component of the complement of
the punctures and double points has negative Euler characteristic.

The spaces $\Mc_{g,n}$ form a modular operad, see Getzler-Kapranov
\cite{modular}, that is, two kinds of operations are defined: attaching
two curves at punctures: $\Mc_{g,n} \times \Mc_{g',n'} \to \Mc_{g+g',
n+n'-2}$ and gluing two punctures together on a single curve $\Mc_{g,n}
\to \Mc_{g+1,n-2}$. Furthermore, there is the action of the symmetric
group $S_n$ on $\Mc_{g,n}$ which reorders the punctures.

Let $F_p = F_p(g,n) \subset \Mc_{g,n}$ be the closed subspace   (substack, in
fact) of dimension $p$ formed by stable curves with at least $\dim_\nc\, (
\Mc_{g,n}) - p = 3g-3+n -p$ double points.  We obtain an ascending filtration
of the moduli space $\Mc_{g,n}$:
\begin{equation*}
\F{-1}\, =\, \emptyset \subset \F{0} \subset \cdots \subset \F{3g-3 +
n}\,=\, \Mc_{g,n}.
\end{equation*}
As in the genus zero case, this filtration behaves nicely with respect to the
modular operad operations:
\[
F_{p}(g,n) \times F_{p'} (g',n') \to F_{p+p'} (g+g', n+n'-2)
\]
corresponding to attaching two curves at punctures and
\[
F_p (g,n) \to F_{p} (g+1, n-2)
\]
corresponding to glueing together two punctures on a single curve.

Irreducible components (strata) $S_G$ of $F_p$ are indexed by {\it stable\/}
labeled $n$-graphs $G$ with $3g - 3 + n  - p + 1$ vertices and with the
invariant $g(G)$, defined below, equal to the genus $g$. {\it Stable\/} refers
to graphs of the following kind. Each graph is connected, has a root vertex and
$n$ enumerated exterior edges, edges which are coincident with only one vertex
of the graph. Each vertex $v$ of the graph is labeled by a nonegative integer
$g(v)$, called the genus of a vertex. The stability condition means that any
vertex $v$ labeled by $g(v)=1$ should be coincident with at least one edge
(i.e., be of at least valence one) and each vertex $v$ with $g(v) = 0$ should
be of at least valence three. The invariant $g(G)$ is given by the formula
$g(G) = b_1(G) + \sum_v g(v)$, where $b_1(G)$ is the first Betti number of the
graph. Each component stratum $S_G$ is a quotient of the product of
uncompactified  moduli
spaces (via the modular operad structure), the combinatorics of which are
neatly encoded in the graph:
\[
S_G = \left( \prod_{v \in G} \M_{g(v), n(v)} \right)/ \aut
\]
where the product is over all vertices of $G$ and $n(v)$ is the valence of the
vertex $v$ and $\aut$ is the automorphism group of a graph $G$ (a bijection on
vertices and edges, preserving the exterior edges, the labels of vertices and
the incidence relation).

Thus, we get a modular operad in the category of filtered varieties (stacks,
in fact) and therefore, applying the spectral sequence functor, we obtain
a modular operad of spectral sequences, see Section~\ref{sec:filter} for
related formalism. Its $(g,n)$-component can be described, as usual, in the
following way
\[
E^1_{p,q} = H_{p+q} (F_p, F_{p-1}, \nc) = H^{p-q} (F_p \setminus F_{p-1},
\nc),
\]
due to Poincar\'{e}-Lefschetz duality, with the differential $d^1:
H^{p-q} (F_p \setminus F_{p-1}) \to H^{p-q-1} (F_{p-1} \setminus
F_{p-2})$ being the Poincar\'{e} residue, and
\[
\bigoplus_{p+q = k} E^\infty_{p,q}
= H_{k} (\Mc_{g,n}, \nc).
\]

\begin{tth}
\begin{enumerate}
\item The $E^1$ term of the spectral sequence is naturally isomorphic to
the Feynman transform, see \cite{modular}, of the modular co-operad
$H^\bullet
(\M_{g,n})$. Namely, $E_{p,q}^1 = 0$, unless $-p \le q \le p \le 3g-3+ n$,
when
\[
E_{p,q}^1 = \bigoplus_G \; \left( \bigoplus_{\sum k(v) = p - q} \quad
\bigotimes_{v \in G} H^{k(v)} (\M_{g(v), n(v)}) \right) ^{\aut},
\]
the first summation running over all stable labeled $n$-graphs $G$ with $g(G) =
g$ and $3g -3 + n -p +1$ vertices, the second over all functions $k(v) \in
\nz$ of vertices $ v \in G$ summing up to $p-q$. The differential is induced by
contracting internal edges in $G$, which corresponds to forming new double
points on a curve, and taking the Poincar\'{e} residue.
\item
When $g = 0$, the graphs $G$ have no loops and are labeled with zeroes,
i.e., just trees where all vertices have valence $\ge 3$. Then the $E^1$
term is nothing but the cobar construction of the co-operad $H^\bullet
(\M_{0,n})$. Moreover, $E_{p,q}^1 = 0$ unless $0 \le q \le p \le n-3$, and
$E^2_{p,q} = E^\infty_{p,q} = 0$, except
\[
E^2_{p,p} = E^\infty_{p,p} = H_{2p} (\Mc_{0,n}).
\]
\end{enumerate}
\end{tth}

\begin{rem}
The second part of the theorem, as well as of Theorem~\ref{dual} below, is
proved independently by Getzler. Just after we proved it, we came across his
two-day old preprint \cite{g:new} dedicated to this kind of duality between the
operads $H_\bullet (\M_{0,n})$ and $H_\bullet(\Mc_{0,n})$.

In contrast to Getzler's proof, our proof below of the statement about
$E^2$ being concentrated on the diagonal uses the operad structure and
known results on the cohomology of $\Mc_{0,n}$.

Getzler \cite{g:new}, on the other hand, uses the fact that the mixed Hodge
structure on the
cohomology of $\M_{0,n}$ is pure and that the operad $H_\bullet (\Mc_{0,n})$ is
Koszul. The Hodge structure on
the cohomology of $\M_{g,n}$ is no longer pure for higher genera, but Mumford's
conjecture for the stable cohomology implies that the Hodge structure on the
stable cohomology of $\M_{g,n}$ should be pure.
\end{rem}

\begin{quest}
For $g > 0$, describe the locus of the $E^\infty$ term on the $(p, q)$ plane.
Is there a stable version of this theorem, that is, for large $g > 0$, where
the sequence degenerates at $E^2$ and the $E^2$ term is located on the diagonal
$p = q$? It would be interesting to find the proper notion of Koszulness for
the stable homology of the modular operad $H^{\operatorname{st}}_\bullet
(\Mc_{g,n})$.
\end{quest}

\begin{pf}
The first part of the theorem is obvious after the description of the strata
$S_G$ above. In the second part, to show that $E_{p,q}^1 = 0$ for $q < 0$
notice that $H^{k(v)} (\M_{0,n(v)} ) = 0 $ for $k(v) > n (v) - 3$.
Denote by $\operatorname{ed} (G)$ the number of edges of the graph,
including the $n$ exterior edges  and by  $v(G)$ the number of vertices.
Adding the above inequalities together, we get $E_{p,q}^1 = 0$
for $p - q > \sum_v
(n(v) - 3) = 2 \operatorname{ed} (G) - n - 3 v(G) = 2 (n + v(G) -1) -n -3 v(G)
= n- 2 - v(G) = p$, that is, for $q < 0$.

The degeneration $E^2_{p,q} = E^\infty_{p,q}$ of the spectral sequence
for $g=0$ follows from the purity of the Hodge structure on $E^1$,
which is the sum of tensor products of the cohomologies of $\M_{0,n}$,
where the Hodge structure is pure, see \cite{bg:1}.

We will show the vanishing $E^2_{p,q} = 0 $ for $g= 0 $ and $p \ne q$
by using Keel's description \cite{keel} of the homology of
$\Mc_{0,n}$, the operad structure and induction on $n$. $\Mc_{0,3}$ is
a point and the statement is trivial. Assuming it is proved for $k \le
n$, let us prove it for $k = n+1$.  Keel's description says that
$H_\bullet (\Mc_{0,n})$ is generated as an intersection algebra by
$H_2$ so $H^1$ is zero.  Thus the term $E^2_{1,0}$ must be zero and
hence the statement is true for $\Mc_{0,4}$.  Now let $E^r_{p,q} (n)$
refer to the spectral sequence for $\Mc_{0,n}$.  Except for the
fundamental class of $\Mc_{0,n+1}$, which is in $E^2_{n-2, n-2}
(n+1)$, the rest of the terms $E^2_{p,q} (n+1)$ are the operad
compositions of $E^2_{p_i, q_i} (n_i)$'s for $n_i \le n$, where we
know $p_i = q_i$ by the induction assumption.
\end{pf}

The situation with the ``Koszul dual'' spectral sequence, the one
which converges to $H^{\bullet} (\M_{g,n}, \nc)$ and whose $E_1$ term
is formed by the cohomology of closed
strata is somewhat better. A similar theorem holds and moreover, the
spectral sequence degenerates at $E_2$ even for $g > 0$, due to the
purity of the Hodge structure on $H^{\bullet} (\Mc_{g,n}, \nc)$, see
Deligne \cite{del:hodge}.

Let $X = \Mc_{g,n}$, $U = \M_{g,n}$ and $D= X \setminus U$. Following
Deligne \cite{del:hodge}, consider the double complex $\Omega^{\bullet,
\bullet}_X (\log D)$, the smooth global $(p,q)$-forms on $X$ with at
most logarithmic singularities $\partial f/f$ along $D$, where $f=0$ is
a local equation of $D$. Its hyper(=total)cohomology is equal to $H^{p+q} (U,
\nc)$. Let $W_m$ be the subcomplex of the total complex generated by
products $\alpha \wedge
\partial f_{i_1}/ f_{i_1} \wedge \dots \wedge \partial f_{i_s}/
f_{i_s}$ for $s \le m$ and smooth $\alpha$, where each $f_i$ is a local
equation of an irreducible component of $D$. The $W^{m} = W_{-m}$'s
for $m \le 0$ define a decreasing filtration of the logarithmic double
complex. Consider the spectral sequence $E_r$ associated with this
filtration. With respect to this filtration, we have $E_1^{p,q} =
H^{2p+q} (\widetilde F_{3g- 3 +n + p})$, where $\widetilde F_s$ is the
disjoint union of irreducible components of $F_s$. The
sequence degenerates: $E_2 = E_\infty = H^{\bullet} (\M_{g,n}, \nc)$,
see \cite{del:hodge}.

In addition, for $g=0$, it is known from Arnold's explicit
description of cohomology classes of configuration spaces that all
nonzero cohomology classes in $H^{m} (\M_{0,n}, \nc)$  are
represented by $(m,0)$-forms with exactly $m$ logarithmic
singularities. Thus, $E^{p,q}_2 = E^{p,q}_\infty$ vanishes unless $p =
-m$ and $ 2p +q = 0$. It follows that $H^m (\M_{0,n}, \nc) =
E^{-m,2m}_2$, where the Hodge structure is pure of weight $2m$, see
\cite{del:hodge}.

  From these descriptions, we have the following theorem.

\begin{tth}
\label{dual}
\begin{enumerate}
\item The $E_1$ term of the spectral sequence is dual to
the Feynman transform of the modular co-operad $H^\bullet
(\Mc_{g,n})$.  Namely, $E^{p,q}_1 = 0$, unless $-3g +3 -n \le p \le
0$ and $-2p \le q \le 6g -6 + 2n$, when
\[
E^{p,q}_1 = \bigoplus_G \; \left( \bigoplus_{\sum k(v) = 2p + q} \quad
\bigotimes_{v \in G} H^{k(v)} (\Mc_{g(v), n(v)}) \right) ^{\aut},
\]
the first summation running over all stable labeled $n$-graphs $G$ with $g(G) =
g$ and $- p +1$ vertices, the second over all functions $k(v) \in \nz$ of
vertices $ v \in G$ summing up to $2p +q$. The differential $d_1: E^{p,q}_1 \to
E^{p+1,q}_1$ is induced by contracting internal edges in $G$, which corresponds
to forming new double points on a curve. The $E_2$ term is equal to $E_\infty =
H^{\bullet} (\M_{g,n}, \nc)$.
\item
When $g = 0$, the graphs $G$ have no loops and are labeled with
zeroes, i.e., just trees where all vertices have valence $\ge 3$. Then
the $E_1$ term is nothing but the dual of the cobar construction of
the co-operad $H^\bullet (\Mc_{0,n})$.
Moreover, $E_2^{p,q} =
E_\infty^{p,q} = 0$, except
\[
E_2^{p,-2p} = E_\infty^{p,-2p} = H^{-p} (\M_{0,n}), \quad 0 \le -p \le n-3.
\]
\end{enumerate}
\end{tth}

\section{Filtered Topological Gravity and \CHA s}

\begin{df}
Let $\Mc$ be the operad of the moduli space of stable curves of genus zero. A
{\sl topological gravity} consists of a differential graded complex vector
space $Q:V_g\,\to\, V_{g-1}$ and a collection $\{\omega_n\}$ such that
\begin{enumerate}
\item $\omega_n = \sum_{r=0}^{2n-4}\omega_n^r$ where $\omega_n^r$ belongs to
$\Omega(\Mc(n))\otimes\End{V}(n)$ with bidegree $(r,-r)$.
\item $d\omega_n = Q\omega_n$
\item $\sigma^*\omega_n = \omega_n\circ \sigma$ for all $\sigma$ in $S_n$
where $\sigma^*$ denotes the pullback of the action of $S_n$ on $\Mc(n)$ and
$\sigma$ on the right hand side indicates the action of $S_n$ on $V^{\otimes
n}$ and
\item $\gamma_i^*\omega_{n+n'-1} = \omega_n\oo_i\omega_n'$
for all $i=1,\ldots, n$ and all $n,n'$ where $\oo_i$ denotes composition in
the endomorphism operad and $\gamma_i^*$ is the pullback of the composition
map $\gamma_i:\Mc(n)\times\Mc(n')\,\to\,\Mc(n+n'-1)$.
\end{enumerate}

A {\sl filtered topological gravity} is a topological gravity in which the
differential forms $\omega_n^p$ vanish for all $p > \dim_\nc \Mc(n) = n-2$.
\end{df}

Since $\Mc$ is an operad of topological spaces, $(C_\bullet(\Mc),\partial)$,
the singular chains on $\Mc$,  inherits the structure of an operad. If $(V,Q)$
is a topological gravity with differential forms $\{\,\omega_n\,\}$ then $V$
is an algebra over $C_\bullet(\Mc)$ where the morphism
$C_\bullet(\Mc(n))\,\to\,\End{V}{n}$ is given by $c\,\mapsto\, \int_c
\omega_n$. The axioms of a topological gravity are formulated precisely so as
to insure that this map is a morphism of operads. This definition of
topological gravity was motivated by its origins in the physics
literature \cite{W}.

Let $(V,Q)$ be a filtered topological gravity with differential forms
$\{\omega_n\}$. Integration of $\{\omega_n\}$ makes $(V,Q)$  into a
$C_\bullet(\Mc)$ algebra with a morphism $\mu: C_\bullet(\Mc)\,\to\,
\End{V}$ such that $\mu(c) = 0$ for all $p$-chains $c$ on $\Mc(n)$ with $p >
\dim_\nc \Mc(n) = n-2$.

\begin{tth}
Let $(V,Q)$ be a filtered topological gravity, let $V$ be filtered   by
the degree in the complex and let $F = \{\,F(n)\,\}$ be the
filtration of $\Mc$ arising from the canonical statification of $\Mc$, then
$V$ is a $C_\bullet(\Mc)$-algebra preserving the filtration.
\end{tth}
\begin{pf}
The filtration of $\Mc$ by $F$ induces a filtration on its singular chains,
\ie\   $\cdots \subseteq C_{p+q}(F_p(n))\subseteq
C_{p+q}(F_{p+1}(n))\subseteq \cdots.$ Similarly, $\End{V}$ is filtered by
$\cdots\subseteq F_p\,\End{V}^{p+q}(n) \subseteq F_p\,\End{V}^{p+q}(n)\subseteq
\cdots$. However, with the given choice of filtration degree, we see that
\begin{equation*}
F_p\,\End{V}(n)^{p+q}(n) = \cases 0, &\text{if $q\geq 1$}\\
\End{V}^{p+q}(n), &\text{if $q\leq 0$.}
\endcases
\end{equation*}
Therefore, the morphism $\mu:C_\bullet(\Mc)\,\to\,\End{V}$ is an algebra
preserving the filtrations if and only if $\mu(c) = 0$ when acting upon
elements in  $C_{p+q}(F_p(n))$ for all $q\geq 1$. This is precisely the case
when $(V,Q)$ is a filtered topological gravity.
\end{pf}

\begin{crl}
A filtered topological gravity is a \CHA.
\end{crl}
\begin{pf}
Let the $E^r$ terms in the spectral sequences associated to $C_\bullet(\Mc)$
and $\End{V}$ be denoted by $E^r$ and $E^r[V]$, respectively. Since $V$ is a
filtered $C_\bullet(\Mc)$-algebra, $E^1[V]$ is an algebra over the operad
$E^1$. However, since $E^1[V] \simeq \End{V}$ and $E^1[V]$ contains a
suboperad $\D^1_{0} = \oplus_{p} H_{p}(F_p(n),F_{p-1}(n))$ which is precisely
the $C_\infty$ operad, $V$ is a \CHA.
\end{pf}

Notice that we only needed the fact that $\mu$ vanished when acting upon
chains that are greater than half the dimension of the corresponding moduli
space, thereby proving Theorem \ref{th:goal}.

\section{Filtrations of Operads and Algebras}
\label{sec:filter}
We will now study properties of operads with a filtration and algebras which
respect this filtration. We shall see that a filtered operad makes each
term in its associated spectral sequence into an operad. Furthermore, there
are natural suboperads in each term of the spectral sequence, one of which
is the \CHA\  in the $E^1$ term of the spectral sequence to the moduli space
of stable curves. All of the results in this section can be naturally
extended to the moduli spaces of higher genus curves in a straightforward
way using the formalism of modular operads.
\begin{df}
Let $\OO = \{\,\OO(n)\,\}$ be an operad of complexes with differentials
$\partial : \OO_p(n)\, \to \, \OO_{p-1}(n)$.  Let $F_p = \{\,F_p(n)
\,\}_{n\geq 1}$ be a filtration of $\OO(n)$ as complexes such that for all
$i=1,\ldots n$, the composition maps $\oo_i$ take $F_{p,q}(n)\otimes
F_{p',q'}(n')$ to $F_{p+p',q+q'}(n+n'-1)$ for all $p, p'$,
$q,q'$,  and the filtration degree $q$ part of $F_p(n)$) is stable
under the action of the permutation group which commutes with the
differential. The collection $F$ is said to be a {\sl
filtration of the operad $\OO$.}
\end{df}
It is clear that a filtered operad can be defined in the category of
topological spaces as well and that such a filtration induces a filtration on
the associated operad of singular chains.

\begin{prop}
Let $\OO$ be an operad with filtration $F$ and let $E^r = \{\, E^r(n)\,\}$ be
the $E^r$ term in the associated spectral sequence then for all $r$, $E^r$
inherits the structure of an operad of complexes with a
differential $\partial^r:E^r_{p,q}(n)\,\to\,E^r_{p-r,q+r-1}(n)$ and
composition maps satisfying
\begin{equation*}
\oo_i:E^r_{p,q}(n) \otimes E^r_{p',q'}(n)\,\to\, E^r_{p+p',q+q'}(n + n' -1)
\end{equation*}
for all $i=1,\ldots, n$. In particular, $E^0_{p,q}(n) =
F_{p,p+q}(n)/F_{p-1,p+q}(n)$ and $E^1_{p,q}(n) = H_{p+q}(F_p(n),
F_{p-1}(n))$.
\end{prop}
\begin{pf}
Let $Z_{p,q}^r(n) = \{\,x\in F_{p,q}(n)\,|\,\partial x\in F_{p-r,q+r-1}(n)
\,\}$ then we can write (see, for example, \cite{lang})
\begin{equation*}
E^r_{p,q}(n) = \frac{Z^r_{p,q}(n)}{\partial Z^{r-1}_{p+(r-1), q- (r-1) +
1}(n) + Z^{r-1}_{p-1,q+1}(n)}
\end{equation*}
with the differential $\partial^r:E^r_{p,q}(n)\,\to\,E^r_{p-r,q+r-1}(n)$
induced from $\partial$. Keeping track of the composition maps, a computation
shows that the numerators assemble into an operad and the denominators into
an operad ideal so their quotient is an operad.
\end{pf}

Given a filtered operad, there is a natural collection of suboperads of each
$E^r$ term as shown by the following result.

\begin{prop}
Let $\OO$ be an operad with filtration $F$. For each $r$ and $k$, there
exists a suboperad of $E^r$, $\D^r_{k} = \{\,\D^r_{k}(n)\,\}_{n\geq 1}$
such that
\begin{equation*}
\D^r_{k}(n) = \bigoplus_{(r-1)p+rq = k(n-1)} E^r_{p,q}(n).
\end{equation*}
In particular, $\D^1_{0} = \{\,\D^1_{0}(n)\,\}$ is a suboperad of $E^1$ with
$n$th component $\D^1_{0}(n) = \bigoplus_p H_{p}(F_p(n),F_{p-1}(n))$.
\end{prop}

\begin{pf}
Each $E^r(n)$ is the union of subcomplexes $K^r_s(n) = \oplus_{(r-1) p + r q =
s} E^r_{p,q}(n)$ since the differential maps $\partial^r:E_{p,q}(n)\, \to\,
E_{p-r,q+r-1}(n)$. The composition maps take $\oo_i: K^r_s(n)\otimes
K^r_{s'}(n') \,\to\, K^r_{s+s'}(n+n'-1)$. Therefore, ${\displaystyle
\D^r_{k}(n)
\linebreak[1] = \linebreak[0] \bigoplus_{(r-1)p + rq = k(n-1)} K^r_{k(n-1)}(n)}
$
assemble to form an operad of complexes.
\end{pf}

In the case of $\Mc$ with its canonical filtration, the suboperads $\D^1_{k}$
are all trivial for dimesional reasons
except when $k$ is $0$.

We now introduce a filtration on the endomorphism operad of $V$ by endowing
$V$ with an additional grading. Let $V$ be a complex with differential
$Q:V^g\,\to\,V^{g-1}$ and a second grading called   the {\sl filtration
degree} denoted by $V_p$ such that $V$ is now a filtered complex. The
endomorphism operad $\End{V}$ inherits a natural filtration $F_p\,\End{V} =
\{\,F_p\,\End{V}(n)\,\}$ where $F_p\,\End{V}(n)$ is the space of maps in
$\End{V}(n)$ with filtration degree less than or equal to $p$. In this case,
the spectral sequence associated to this filtration makes each $E^r$ term
into an operad as well as the suboperads indicated above. Perhaps the
simplest filtration of $V$ is when the filtration degree is exactly the
degree in the complex. In this case, the associated spectral sequence
degenerates at the $E^2$ term, where $E^0 \simeq \End{V}$ with a zero
differential, $E^1 \simeq\End{V}$ with the differential $Q$, and $E^2 \simeq
H(\End{V})$.

\begin{df}
Let $V$ be a vector space graded by a filtration degree as above. We say that
an
$\OO$-algebra, $V$, is a {\sl filtered $\OO$-algebra} if the morphism of
operads $\OO\,\to\, \End{V}$ preserves the filtrations.
\end{df}

If $\OO$ is an operad with filtration $F$ and $V$ a filtered $\OO$-algebra,
then the $E^r$ term of $\End{V}$ is an algebra over the $E^r$ term of
$\OO$ for all $r\geq 0$.

In the case of the operad $\Mc$ with its canonical filtration where $V$ is a
filtered $C_{\bullet}(\Mc)$-algebra with  the filtration degree on $V$
equal to  the degree of the complex, then $V$ is an algebra (with zero
differential) over the operad $E^0$, $V$ is an algebra (with $Q$
differential), over the operad $E^1$, and $H(V)$ is an algebra over $E^2$
(which is isomorphic to the operad $H_\bullet(\Mc)$) where we have used the
canonical morphism $H(\End{V})\,\to\,\End{H(V)}$ in the last step.

\begin{ack}
We would like to thank F.~Cohen, M.~Kontsevich, R.~Hain and
M.~Schlessinger for helpful discussions on the cohomology of moduli
spaces. We thank Gregg Zuckerman for his useful comments on Table~1. The third
author is very grateful to the hospitality of the
Weizmann Institute of Science, where part of the work has been done.
\end{ack}

\bibliographystyle{amsplain}

\makeatletter \renewcommand{\@biblabel}[1]{\hfill#1.}\makeatother

\end{document}